\documentclass[draft]{amsart}
\usepackage{amsmath, amssymb, epsfig, graphicx}

\renewcommand{\P}{\mathcal P}

\renewcommand{\leq}{\leqslant}
\renewcommand{\geq}{\geqslant}
\newcommand{\supp}{{\sf Supp}\,}

\newcommand{\gen}[1]{\bigl<#1\bigr>}

\newcommand{\sym}[1]{{\sf S}{(#1)}}
\newcommand{\trans}[1]{{\sf T}(#1)}
\newcommand{\transs}[2]{{\sf T}(#1,#2)}
\newcommand{\alt}[1]{{{\sf A}{(#1)}}}
\renewcommand{\wr}{\,{\sf wr}\,}
\newcommand{\id}{{\rm id}}
\newcommand{\rank}{\mbox{\emph{\mbox{\rm rank}\,}}}
\newcommand{\End}[1]{\mbox{End}(#1)}
\newcommand{\F}{\mathbb F}

\newtheorem{thm}{Theorem}[section]
\newtheorem{theorem}{Theorem}[section]
\newtheorem{lem}[thm]{Lemma}
\newtheorem{lemma}[thm]{Lemma}

\newtheorem{exam}[thm]{Example}
\newtheorem{defn}[thm]{Definition}

\title[Rank of the Endomorphism Monoid of a Partition]
{The Rank of the Endomorphism Monoid of a Partition}
\author{Jo\~ao Ara\'ujo \and Csaba Schneider}
\address[Schneider]{Informatics Research Laboratory\\
Computer and Automation Research Institute\\
1518 Budapest Pf.\ 63\\
Hungary\hfill\break
Centro de \'Algebra da Universidade de Lisboa\\
Av.\ Prof.\ Gama Pinto 2\\ 1649-003 Lisboa\\ Portugal\hfill\break
Email: csaba.schneider@gmail.com\hfill\break
WWW: www.sztaki.hu/$\sim$schneider}
\address[Ara\'ujo]{Universidade Aberta\\ R. Escola Polit\'{e}cnica 147\\
1269-001 Lisboa\\ Portugal\hfill\break
Centro de \'{A}lgebra\\ Universidade de Lisboa\\
1649-003 Lisboa\\ Portugal\hfill\break
Email: jaraujo@lmc.fc.ul.pt}



\date{\today}

\begin{document}
\begin{abstract}
The rank of a semigroup is the cardinality of a smallest generating set. In
this paper we compute the rank of the endomorphism monoid of a non-trivial uniform partition of a finite set, that is, the semigroup of those transformations of a
finite set that leave a non-trivial uniform partition invariant. That involves proving that the rank of a wreath product of two symmetric groups is two and then  use the fact that the endomorphism monoid of a partition  is isomorphic to a wreath product of two full transformation
semigroups. The calculation of the rank of these semigroups solves an open question.  
\end{abstract}
\maketitle

\section{Introduction}

If $S$ is a semigroup and $U$ is a
subset of $S$ then we say that {\em $U$ generates $S$} if every element of $S$ is
expressible as a word in the elements of $U$. We use the convention that
the empty word represents the
identity element. The {\em rank} of a
semigroup $S$, denoted by $\rank S$,
is the minimum among the cardinalities of its generating sets.
It is well-known that a finite full transformation semigroup
has rank~3, while a finite full partial transformation semigroup
has rank~4 (see~\cite[Exercises~1.9.7 and 1.9.13]{howie}). 
Similar results were proved for many different classes of 
transformation semigroups (such as total, partial, partial one-to-one, order preserving) and their ideals; see \cite{1,8,10,15,16}. 

Some generalizations of the notion of  rank, for instance the {\em idempotent rank} and the {\em nilpotent rank}, also attracted a great deal of attention (see \cite{2,6,14}, among others).

Finally, in recent years, the notion of relative rank  has been subjected to extensive research (see for example \cite{0,mitchell1,8.5,9,mitchell2,11}). Relative rank is a useful notion when dealing with finite semigroups (see Lemma~\ref{l1}), and it is crucial when dealing with uncountable semigroups. In fact, in such semigroups, the notion of rank is not very informative as  the rank and the 
cardinality coincide.
To a large extent, this line of research was prompted by two old papers by Sierpi\'nski  \cite{si} and by Banach \cite{sb-1935} (see also \cite{ams}).

In this paper we deal with the  endomorphism monoid of a non-trivial uniform partition. We prove that the rank of such semigroup is 4, thus settling a problem posed in~\cite{huisheng}. We also calculate
the ranks of some related transformation semigroups.


If $X$ is a finite set then the set of transformations on $X$ form the {\em
  full transformation monoid on $X$} and is denoted by $\trans X$. We assume that
  transformations act on the right; that is, if $x\in X$ and $f\in \trans X$, then
$xf$ denotes the image of $x$ under $f$.
Let $\P$ be a partition of $X$; that is,
$\mathcal P=\{P_1,\ldots,P_m\}$ where $P_1,\ldots,P_m\subseteq X$,
$P_i\cap P_j=\emptyset$ whenever $i\neq j$, and $X=P_1\cup\cdots\cup
  P_m$. The equivalence relation that corresponds to a partition $\P$ is
  denoted $\sim_\P$. The elements of $\trans X$ that preserve the partition $\P$ form
 a semigroup and is denoted by $\transs X \P$. Using symbols,
\[
\begin{array}{rcl}
  \transs X\P &=&  \{f\in \trans X \mid (\forall{P_i\in \P})(\exists{P_j\in\P})
\ P_if\subseteq P_j\} \\
   & = & \{f\in \trans X \mid \mbox{if }x\sim_\P y \mbox{ then }
xf\sim_\P yf\}.
\end{array}
\]
The partition $\P$ is called {\em uniform} if $|P_i|=|P_j|$ for all $i,\ j\in\{1,\ldots,m\}$.
The main result of Huisheng's paper~\cite{huisheng} is that for a uniform
partition $\P$ the rank of $\transs X\P$ is at most~6,
or smaller in some degenerate cases.
Huisheng's proof relied on the observation that there is a strong relationship
between $\rank \transs X\P$ and the rank of the group $G$ of invertible
transformations in $\transs X\P$. He proved that the rank of $G$ is at
most 4 (or smaller in some degenerate cases).
In the present paper, we are able to show, for $|X|\geq 3$, that $\rank G=2$
(see Theorem~\ref{wrth}). In order to facilitate the
 proof of our results, we
use the concept of wreath products of transformation semigroups; see
Section~\ref{wreath}.

We can also consider the rank of some further interesting semigroups related
to a partition $\P$ of a set $X$. Let us define
\[
\begin{array}{rcl}
  \Sigma(X,\P)    & = & \{f\in \trans X\ |\  x\sim_\P y \mbox{ if and only if }
xf\sim_\P yf\},
\mbox{ and} \\
  \Gamma(X,\P)    & = & \{f\in \transs X\P \mid (\forall{P_i,\ P_j\in \P})\ \mbox{ either
  } P_if\cap P_j\neq
 \emptyset \mbox{ or }P_if=P_j \}. \\
\end{array}
\]
It is routine to check that $\Sigma(X,\P),\ \Gamma(X,\P)\leq \transs X\P$.

The following theorem improves Huisheng's
result and settles the problem of determining the rank of $\transs X\P$. It
also gives the rank of $\Sigma(X,\P)$ and $\Gamma(X,\P)$. A partition of $X$ is said to be trivial if it has $1$ or $|X|$ parts (that is, the trivial partitions are the identity and the universal partitions). 

\begin{theorem}\label{main}
If $X$ is a finite set such that $|X|\geq 3$, and $\P$
is a non-trivial uniform partition of $X$, then $\rank \transs X\P=4$ and
$\rank\Sigma(X,\P)=\rank \Gamma(X,\P)=3$.
\end{theorem}

\section{Wreath products of transformation semigroups}\label{wreath}

We define the wreath product of two transformation semigroups
following~\cite{straubing} (see also~\cite[Chapter~10]{mel}). The material of this section is
well-known, however, we felt that it was necessary to present it
in order to make the paper
 self-contained and also to set our system of notation.
Let $S$ and $R$ be two transformation semigroups
acting on the sets $Y$ and $Z$, respectively, and let $B$ denote the set
of functions $f:Z\rightarrow S$. The underlying set of the wreath product
$S\wr R$ is the Cartesian product $B\times R$ and to each element $(f,r)\in
B\times R$, we assign a transformation of the set
$Y\times Z$ defined by
\begin{equation}\label{wract}
(y,z)(f,r)=(yf(z),zr)\quad\mbox{for all}\quad y\in Y,\ z\in Z,\ f\in B,\ r\in R.
\end{equation}
It is easy to see that this assignment is injective and that the image of this
assignment is closed under composition; that is, the image is  a subsemigroup of $\trans{Y\times Z}$.

From now on we assume that the sets $Y$ and $Z$ are finite and that $S$ and
$R$ are monoids, that is, they contain the identity transformation $\id$.
Let us assume without loss of generality that $Y=\{1,\ldots,n\}$ and
$Z=\{1,\ldots,m\}$.
In this
case, a function $f:Z\rightarrow S$ can be represented as the $m$-tuple
$(s_1,\ldots,s_m)$ where $f(i)=s_i$ for all $i\in\{1,\ldots,m\}$. This defines
a bijection between $B\times R$ and $S^m\times R$, and, in turn, between
$S\wr R$ and $S^m\times R$. Therefore we may view an element of $S\wr R$ as a
pair $((s_1,\ldots,s_m),r)$ where $s_1,\ldots,s_m\in S$ and $r\in R$. The
element $((s_1,\ldots,s_m),r)$ will more briefly be denoted by $(s_1,\ldots,s_m)r$.
Setting $s=(s_1,\ldots,s_m)$ the same element can be written as $sr$.
By~\eqref{wract},
the action of $(s_1,\ldots,s_m)r$ on a pair $(y,z)\in Y\times Z$ is given by
\begin{equation}\label{wract2}
(y,z)(s_1,\ldots,s_m)r=(ys_z,zr).
\end{equation}

Let $S\leq \trans Y$ and $R\leq\trans Z$ be transformation monoids as above.
Then $S$ and $R$ can naturally be embedded into
the wreath product $S\wr R$.
Indeed, we may consider the following submonoids of $S\wr R$:
\begin{equation}\label{si}
\overline S_i=\{(\id,\ldots,\id,\stackrel{\mbox{\footnotesize $i$-th
    component}}{s},\id\ldots,\id)\id \in S\wr R\ |\ s\in S\},
\end{equation}
and
\begin{equation}\label{r}
\overline R=\{(\id,\ldots,\id)r \in S\wr R\ |\
r\in R\}.
\end{equation}
Set, for $i=\{1,\ldots,m\}$, $\overline
Y_i=\{(y,i)\ |\ y\in Y\}$. Easy computation shows that the elements of
$\overline S_i$ leave the
set $\overline Y_i$ invariant, and $\overline S_i$, considered as a
transformation monoid on $\overline Y_i$,
is isomorphic to the transformation monoid $S$.
Let us further define $\overline Z=\{\overline Y_1,\ldots,\overline Y_m\}$. Then
$\overline R$ leaves $\overline Z$ invariant, and $\overline R$, considered as
a transformation monoid on $\overline Z$, is isomorphic to $R$.
Since $\overline S_1\times\cdots\times \overline S_m\cong S^m$ and
$\overline R\cong R$, we may consider $S^m$ and $R$ as submonoids of $S\wr R$.

Since $R$ is assumed to be a transformation monoid on the set
$Z=\{1,\ldots,m\}$, we may define a homomorphism $\vartheta:R\rightarrow S^m$:
$$
(r\vartheta)(s_1,\ldots,s_m)=(s_{1r},\ldots,s_{mr}).
$$
We note that the action of $\End{S^m}$ on $S^m$ is a left-action. The
homomorphism $\vartheta$ is useful for expressing the operation in $S\wr
R$. Indeed, let $s_1,\ s_2\in S^m$ and $r_1,\ r_2\in R$. Then, viewing
$s_1r_1$ and $s_2r_2$ as elements of $S\wr R$, easy computation
shows that
\begin{equation}\label{wrop}
(s_1r_1)(s_2r_2)=(s_1 (r_1\vartheta)s_2)(r_1r_2).
\end{equation}
Therefore the wreath product
$S\wr R$ can also be viewed as the semidirect product $S^m\rtimes R$
of $S^m$ and $R$ with respect to the homomorphism
$\vartheta:R\rightarrow\End{S^m}$ (see~\cite[page~186]{mel}).
In particular, if $r$ is an invertible element of $R$ and $s\in S^m$,
then, considering $s$ and $r$ as elements of $S\wr R$, equation~\eqref{wrop} implies that
\begin{equation}\label{conj}
rsr^{-1}=(r\vartheta)s.
\end{equation}
Therefore conjugation by an invertible element of $R$ leaves
$S^m$ invariant, and the conjugation action of $R$ is given
by the homomorphism $\vartheta$. Indeed, from the last displayed line we
obtain that $s^r=r^{-1}sr=(r^{-1}\vartheta)s$. Note
that the conjugation action of the element $r$ is actually given by the
endomorphism $r^{-1}\vartheta$. The reason for this is that we assumed that
endomorphisms act on the left, while the usual definition makes the
conjugation action a right-action.
Further, if $r$ is an invertible element of $R$ and $i,\ j\in\{1,\ldots,m\}$
such that $ir=j$ then~\eqref{conj} implies that
 $\overline S_i^r=(r^{-1}\vartheta) \overline
S_i=\overline S_j$.

For a finite set $X$,
the invertible elements of $\trans X$ form the symmetric group $\sym
X$. If $\P$ is a partition of $X$ then set $\sym{X,\P}=\transs X\P\cap\sym X$.
The next lemma connects the semigroups related to a partition to wreath
products.

\begin{lemma}\label{wreq}
Let $Y=\{1,\ldots,n\}$ and $Z=\{1,\ldots,m\}$, set $X=Y\times Z$, and let $\P$
denote the partition
$\{\{(1,1),\ldots,(n,1)\},\ldots,\{(1,m),\ldots,(n,m)\}\}$.
Then
\begin{enumerate}
\item[(i)] $\transs X\P=\trans Y\wr \trans Z$;
\item[(ii)] $\Sigma(X,\P)=\trans Y\wr\sym Z$;
\item[(iii)] $\Gamma(X,\P)=\sym Y\wr \trans Z$;
\item[(iv)] $\sym{X,\P}=\sym Y\wr\sym Z$.
\end{enumerate}
\end{lemma}
\begin{proof}
Since the proofs of these statements are very similar to each other,
we only verify
assertion~(ii). Let $(y_1,z_1),\ (y_2,z_2)\in Y\times Z$ and assume
that $t_1,\ldots,t_m\in \trans Y$
$s\in\sym Z$ so that
$(t_1,\ldots,t_m)s\in \trans Y\wr\sym Z$. Set
$w=(t_1,\ldots,t_m)s$. We have, for $i=1,\ 2$, that
$$
(y_i,z_i)w=(y_i,z_i)(t_1,\ldots,t_m)s=(y_it_{z_i},z_is).
$$
Thus, if $(y_1,z_1)\sim_\P(y_2,z_2)$, then $z_1=z_2$, and then $z_1s=z_2s$;
hence, in this case, $(y_1,z_1)w\sim_\P(y_2,z_2)w$. Conversely, if
$(y_1,z_1)w\sim_\P(y_2,z_2)w$, then $z_1s=z_2s$, which, using that $s$ is a invertible, gives that $z_1=z_2$;
therefore $(y_1,z_1)\sim_\P(y_2,z_2)$. This shows that $w\in\Sigma(X,\P)$, and
so $\trans Y\wr\sym Z\leq\Sigma(X,\P)$.

Suppose now that $x\in\Sigma(X,\P)$. Then the defining property of
$\Sigma(X,\P)$ implies that $x$ induces a permutation on the  set $\P$. Since
there is a one-to-one correspondence between $\P$ and $Z=\{1,\ldots,m\}$, we
obtain that $x$ induces a permutation on $Z$. Let this
permutation be $s$. For $i=1,\ldots,m$, let us define
a transformation $t_i\in\trans Y$. Let $j\in\{1,\ldots,m\}$ such that
$ix=j$. Then, for all $k\in\{1,\ldots,n\}$ there is some $l_k\in\{1,\ldots,n\}$
such that  $(k,i)x=(l_k,j)$. Let $t_i$ be the transformation that maps $k$
to $l_k$ for all $k\in\{1,\ldots,n\}$. Then routine computation shows that
$x=(t_1,\ldots,t_m)s$, and so $x\in \trans Y\wr\sym Z$. Thus,
$\Sigma(X,\P)\leq \trans
Y\wr\sym Z$, and  so $\trans Y\wr\sym Z=\Sigma(X,\P)$
\end{proof}

\section{Relative rank of semigroups}

Let $U\subseteq S$ be a subset of a semigroup $S$. The
{\em relative rank} of $S$ modulo $U$,
denoted by $\rank{(S:U)}$, is the
minimum among the cardinalities of subsets $V$ of $S$
such that $S=\langle V\cup
U\rangle$.  The relative rank was introduced in~\cite{8.5}.
The next lemma shows that the rank of a transformation semigroup
is related to its relative rank modulo the unit group.

\begin{lem}\label{l1}
Let $S$ be a finite transformation semigroup and let $G$ be the group of units
in $S$. If $U\subseteq S$ such that $\gen U=S$, then $\gen {U\cap G}=G$.
In particular, $\rank S=\rank{(S:G)}+\rank G$.
\end{lem}
\begin{proof}
It suffices to show that  $G\leq \gen {U\cap G}$, and so
suppose that $g\in G$. Since $U$ is a generating set of $S$, we obtain that
$g=u_1u_2\cdots u_r$ with some $u_1,\ldots,u_r\in U$. Since $g$ is invertible,
we obtain that $u_1,\ldots,u_r$ must also be invertible, and so
$u_1,\ldots,u_r\in G$. Thus $g\in\gen{U\cap G}$, and hence the assertion
follows.
\end{proof}

Next we determine the relative ranks of $\transs X\P$, $\Sigma(X,\P)$ and
$\Gamma(X,\P)$ modulo their unit groups.

\begin{lemma}\label{l2}
If $X$ is a finite set and $\P$ is a uniform partition of $X$, then
$$
\rank{(\transs X\P:\sym{X,\P})}=2,
$$
and
$$
\rank{(\Gamma(X,\P):\sym{X,\P})}=\rank{(\Sigma(X,\P):\sym{X,\P})}=1.
$$
\end{lemma}
\begin{proof}
We may suppose without loss of
generality that $X=Y\times Z$ where $Y=\{1,\ldots,n\}$, $Z=\{1,\ldots,m\}$
and $\P$ is the partition
$$
\{\{(1,1),\ldots,(n,1)\},\ldots,\{(1,m),\ldots,(n,m)\}\}.
$$
By Lemma~\ref{wreq},
$\transs X\P=\trans Y\wr\trans Z$, $\Sigma(X,\P)=\trans Y\wr\sym Z$,
$\Gamma(X,\P)=\sym Y\wr\trans Z$, and $\sym{X,\P}=\sym Y\wr\sym Z$, and so it
suffices to show that
\begin{equation}\label{cl1}
\rank{(\trans Y\wr \trans Z:\sym Y\wr \sym Z)}=2
\end{equation}
and
\begin{equation}\label{cl2}
\rank{(\trans Y\wr\sym Z:\sym Y\wr\sym Z)}=\rank{(\sym Y\wr\trans Z:\sym
  Y\wr\sym Z)}=1.
\end{equation}

Let
$\overline \alpha$ denote the transformation in $\trans Y$ such that $1\overline
\alpha=2$
and $i\overline\alpha=i$ for all $i\in\{2,\ldots,n\}$.
Then $\trans Y=\gen{\sym Y\cup \{\overline\alpha\}}$ (see~\cite[Exercise~1.9.7]{howie}). Set $\alpha=(\overline
\alpha,\id,\ldots,\id)\id$ and $\beta=(\id,\ldots,\id)\overline\alpha$.
We claim that
\begin{eqnarray}
\trans Y\wr\sym Z&=&\gen{\sym Y\wr\sym Z\cup\{\alpha\}},\label{eq1}\\
\sym Y\wr\trans Z&=&\gen{\sym Y\wr\sym Z\cup\{\beta\}},\label{eq2}\\
\trans Y\wr\trans Z&=&\gen{\sym Y\wr\sym Z\cup\{\alpha,\beta\}}.\label{eq3}
\end{eqnarray}
For $i\in\{1,\ldots,m\}$, let us define the following submonoids of $\trans Y\wr
\trans Z$:
\begin{eqnarray*}
\overline T_i&=&\{(\id,\ldots,\id,\stackrel{\mbox{\footnotesize $i$-th
    component}}{t},\id,\ldots,\id)\id \in \trans Y\wr \trans Z\ |\ t\in \trans Y\},\\
\overline S_i&=&\{(\id,\ldots,\id,\stackrel{\mbox{\footnotesize $i$-th
    component}}{s},\id\ldots,\id)\id \in \trans Y\wr \trans Z\ |\ s\in \sym Y\},\\
\overline{\trans Z}&=&\{(\id,\ldots,\id)t \in \trans Y\wr \trans Z\ |\ t\in \trans
    Z\},\\
\overline{\sym Z}&=&\{(\id,\ldots,\id)s \in \trans Y\wr \trans Z\ |\ s\in \sym Z\}.\\\end{eqnarray*}
Let us first prove~\eqref{eq1}.
As $\alpha\in\trans Y\wr\sym Z$, we find that $\gen{\sym Y\wr\sym Z\cup\{\alpha\}}\leq \trans Y\wr\sym Z$.
Since
$\trans Y=\gen{\sym Y\cup\{\overline \alpha\}}$, we obtain that $\overline
T_1=\gen{\overline S_1
\cup\{\alpha\}}$,
and hence $\overline T_1\leq \gen{\sym Y\wr\sym Z\cup\{\alpha\}}$. For all $i\in\{1,\ldots,m\}$, there is some
$r\in\sym Z$ such that $1r=i$, and, as discussed before Lemma~\ref{wreq}, we
obtain that $r^{-1}\overline T_1r=\overline T_i$. Therefore $\overline
    T_1,\ldots,\overline T_m\leq \gen{\sym Y\wr\sym Z\cup\{\alpha\}}$. As $\overline{\sym Z}\leq \gen{\sym Y\wr\sym Z\cup\{\alpha\}}$ and
$\trans Y\wr\sym Z=(\overline T_1\times\cdots\times \overline T_m)\rtimes \overline{\sym Z}$, we obtain
that  $\trans Y\wr\sym Z\leq \gen{\sym Y\wr\sym Z\cup\{\alpha\}}$, and so the
    required equality holds.

Now we show that~\eqref{eq2}.
As $\beta\in \sym Y\wr\trans Z$, we have that
$\gen{\sym Y\wr\sym Z\cup\{\beta\}}\leq \sym Y\wr\trans Z$.
As  $\trans Y=\gen{\sym Y\cup\{\overline\alpha\}}$, we obtain that
$\overline{\trans Z}=\gen{\overline{\sym Z}\cup\{\beta\}}$. As
    $\sym Y\wr\trans Z=(\overline S_1\times\cdots\times \overline S_m)\rtimes\overline{\trans Z}$, we
    obtain that $\sym Y\wr\trans Z\leq \gen{\sym Y\wr\sym Z\cup\{\beta\}}$, and so
    the claim is proved.

As $\trans Y\wr\trans Z=(\overline T_1\times\cdots\times \overline
T_m)\rtimes \overline{\trans Z}$,
the arguments in the previous two paragraphs show that
$\trans Y\wr\trans Z\leq \gen{\sym Y\wr\sym Z\cup\{\alpha,\beta\}}$.
As $\alpha,\ \beta\in\trans Y\wr\trans Z$, we obtain~\eqref{eq3}.

As $\sym Y\wr\sym Z$ is a proper submonoid of each of
the monoids $\trans Y\wr\trans Z$, $\trans Y\wr\sym Z$ and $\sym Y\wr\trans Z$,
equation~\eqref{cl2}  must be valid.
In order to show~\eqref{cl1}, it suffices to prove that $\rank{(\trans Y\wr\trans Z:\sym
  Y\wr\sym Z)}>1$. Suppose that $\gamma\in \trans Y\wr\trans Z$ such that
$\gen{\sym Y\wr\sym Z\cup\{\gamma\}}=\trans Y\wr\trans Z$.
Then there are $g,\ g_1,\ldots,g_k\in\sym Y\wr\sym Z$, such that $g\gamma g_1\gamma
\cdots \gamma g_k =\alpha$. Thus
$\gamma g_1\gamma \cdots \gamma g_k =g^{-1}\alpha$ and hence $
\ker \gamma  \subseteq \ker g^{-1}\alpha =
\{((1,1)g,(2,1)g)\}\cup \Delta$, where $\Delta=\{(x,x)\mid x\in X\}$. Since
$\gamma \not\in
\sym Y\wr\sym Z$, we obtain that $\ker \gamma  =
\{((1,1)g,(2,1)g)\}\cup \Delta $.
Similarly, there exist $h,\ h_1,\ldots,h_k\in\sym Y\wr\sym Z$
such that $h\gamma h_1\gamma \cdots \gamma h_k =\beta$,
Thus $\gamma h_1\gamma
\cdots \gamma h_k =h^{-1}\beta$ and hence
\[
\ker \gamma \subseteq \ker h^{-1}\beta = \{((1,1)h,(1,2)h),
\ldots, ((n,1)h,(n,2)h)\}\cup \Delta.
\]
Hence, for some $i\in \{1,\ldots,n\}$, we have $
((1,1)g,(2,1)g)=((i,1)h,(i,2)h)$, that is, $
((1,1),(2,1))gh^{-1}=((i,1),(i,2))$.
This, however, is a contradiction, because
$(1,1)\sim_\P(2,1)$, but $(i,1)\not\sim_\P (i,2)$ and the transformation
$gh^{-1}$ preserves the
equivalence relation~$\sim_\P$.

Therefore we verified that $\gen{\sym Y\wr\sym
  Z\cup\{\gamma\}}$ is a proper submonoid of the wreath product
$\trans Y\wr\trans Z$, which
  shows that equation~\eqref{cl1} must hold.
\end{proof}

\section{The Rank of $\sym{X,\P}$}

In this section we prove the following theorem.

\begin{theorem}\label{wrth}
If $X$ is a finite set such that $|X|\geq 3$ and
$\P$ is a uniform partition of $X$ then $\sym{X,\P}$ is generated by two
elements.
\end{theorem}

In our terminology, the previous theorem gives that $\rank\sym{X,\P}=2$. We
note that the rank of a transitive
permutation group, such as $\sym{X,\P}$, is defined in
permutation group theory as the number of orbits of a point-stabilizer. Thus,
in order to avoid possible confusion,
we decided to state the theorem above without using
the notation $\rank\sym{X,\P}$.

Theorem~\ref{main} will follow from Lemmas~\ref{l1} and~\ref{l2} and
Theorem~\ref{wrth}.

Using the fact that $\sym{X,\P}$ is isomorphic to a wreath product $\sym
Y\wr\sym Z$, it is not difficult to see that $\sym{X,\P}$ is generated by four
elements.  Indeed, consider the subgroups $\overline{\sym Y}_i$ and
$\overline{\sym Z}$ defined in~\eqref{si} and~\eqref{r}. Since
$\sym Y\wr\sym Z=(\overline{\sym Y}_1\times\cdots\times \overline{\sym
  Y}_m)\rtimes
\sym Z$ and
$\sym Z$ is transitive by conjugation on the
subgroups $\overline{\sym Y}_i$, we obtain that $\sym Y\wr\sym
Z=\gen{\overline{\sym Y}_1,\overline{\sym Z}}$. Since
$\overline{\sym Y}_1$ and $\overline{\sym Z}$ are full symmetric groups, they
are generated by two elements, and so we obtain that $\sym Y\wr\sym
Z$ is generated by at most four elements. Essentially this is proved
in~\cite[Theorem~2.6]{huisheng}.

Before proving Theorem~\ref{wrth}, we state two simple lemmas.

Let $G$ be a permutation group acting on the set $\{1,\ldots,n\}$ and let $\F$
be a field. Let $V$ denote the $n$-dimensional vector space over $\F$ with
basis $\{e_1,\ldots,e_n\}$. The group $G$ can be thought of
as a permutation group
on
the set $\{e_1,\ldots,e_n\}$, and this defines an $\F G$-module structure on
$V$ as follows:
$$
e_i g=e_{ig}\quad\mbox{for}\quad i\in\{1,\ldots,n\}\mbox{ and }g\in G.
$$
The module $V$ is called the {\em permutation module} for $G$ over $\F$.

The following lemma is well-known; see, for instance,~\cite[Lemma~5.3.4]{kl}.

\begin{lemma}\label{replem}
If $X=\{1,\ldots,n\}$, then the permutation module for $\sym X$ over a field $\F$ of
characteristic $p$ has precisely two
proper non-trivial submodules:
\begin{eqnarray*}
U_1&=&\left\{(a,a,\ldots,a)\ |\ a\in \F\right\}\quad\mbox{and}\\
U_2&=&\left\{(a_1,\ldots,a_n)\ |\ a_1+\cdots+a_n=0 \right\}.
\end{eqnarray*}
Further, if $p\mid n$ then $U_1\leq U_2$; otherwise $V=U_1\oplus U_2$.
\end{lemma}

Suppose that $G=G_1\times\cdots\times G_k$ where the $G_i$ are
finite groups. For
$I\subseteq\{G_1,\ldots,G_k\}$ the function $\varrho_I:G\rightarrow\prod_{G_i\in
I}G_i$ is the natural projection map. We also write
$\varrho_i$ for $\varrho_{\{G_i\}}$.
A subgroup $X$ of $G$ is said to be a {\em strip} if for
each $i=1,\ldots,k$ either $X\varrho_i=1$ or $X\varrho_i\cong X$. The set
of $G_i$ such that $X\varrho_i\neq 1$ is called the {\em support} of
$X$ and is denoted $\supp X$.
Two strips $X_1$ and $X_2$ are {\em disjoint}
if $\supp X_1\cap\supp X_2=\emptyset$. A strip $X$ is said to be {\em full} if
$X\varrho_i=G_i$ for all $G_i\in\supp X$,
and $X$ is called {\em non-trivial} if
$|\supp X|\geq 2$. A subgroup $K$ of $G$ is said to be {\em subdirect} if
$K\varrho_i=G_i$ for all $i$.

We recall a well-known lemma on finite simple groups which can be
found in~\cite[page~328]{scott}. The proof of the lemma is elementary and does not use
the finite simple group classification.

\begin{lemma}\label{scott1}
Let $M$ be a direct product of finitely many non-abelian, finite simple groups
and let $H$ be a
subdirect subgroup of $M$. Then $H$ is the direct product of pairwise
disjoint full strips of $M$.
\end{lemma}

The wreath product of two transformation semigroups $S$ and $R$ was defined in
Section~\ref{wreath}, and the definition can also be used to construct the
wreath product of two permutation groups $G$ and $H$.
Recall that the wreath product $G\wr H$ is isomorphic to
$G^m\rtimes H$ and
a typical element of $G\wr H$ is denoted
by $(\pi_1,\ldots,\pi_m)\sigma$ where $\pi_i\in G$ and $\sigma\in H$.
Setting $\pi=(\pi_1,\ldots,\pi_m)$, the same element can also be written as
$\pi\sigma$.  The following lemma facilitates the calculations in $G\wr H$.

\begin{lemma}\label{lem}
Let $\pi\sigma,\ \pi_1\sigma_1,\ \pi_2\sigma_2\in G\wr H$ where $G$ and $H$
are as above. Then
\begin{enumerate}
\item[(i)]
  $\pi_1\sigma_1\pi_2\sigma_2=\pi_1(\pi_2)^{\sigma_1^{-1}}\sigma_1\sigma_2$;
\item[(ii)] $(\pi\sigma)^{-1}=(\pi^{-1})^{\sigma}\sigma^{-1}$;
\item[(iii)] $(\pi\sigma)^n=\pi\pi^{\sigma^{-1}}\pi^{\sigma^{-2}}\cdots
\pi^{\sigma^{-n+1}}\sigma^n$.
\end{enumerate}
In particular the projection map $\varrho:G\wr H\rightarrow H$ defined by
$\pi\sigma\mapsto \sigma$ is a homomorphism.
\end{lemma}
\begin{proof}
The assertion that $\varrho$ is a homomorphism follows from~(i).
The rest can be verified using~\eqref{wrop} and~\eqref{conj}.
\end{proof}

Now we are ready to prove Theorem~\ref{wrth}. Permutations of a finite set
will be written as products of disjoint cycles.

\begin{proof}[The proof of Theorem~$\ref{wrth}$]
By Lemma~\ref{wreq}(iv), it suffices to show, for $Y=\{1,\ldots,n\}$ and
$Z=\{1,\ldots,m\}$, that the group $W=\sym Y\wr\sym Z$ is
generated by two elements whenever $nm\geq 3$. Since, for a finite set
$Y$,
the group  $\sym Y$ is generated by two elements, we may assume that $n\geq
2$ and $m\geq 2$.
Let $\alt Y$ denote the group of even permutations of $Y$. Then $\alt Y$ is a
normal subgroup of index two of $\sym Y$.
As $W=\sym Y^m\rtimes \sym Z$,
we may consider the subgroups $\alt Y^m$ and $\sym Y^m$ of $W$ and let $A$ and $S$ denote
these subgroups respectively.

Let us define
\begin{eqnarray*}
x&=&\left\{\begin{array}{ll}
(\id,(1,2),\id,\ldots,\id)(1,2,\ldots,m) & \mbox{if either $n$ or $m$ is
  odd}\\
(\id,(1,2),\id,\ldots,\id)(2,3,\ldots,m) & \mbox{otherwise}\end{array}\right.\\
y&=&((1,2,\ldots,n),\id,\ldots,\id)(1,2).
\end{eqnarray*}
Set $M=\left<x,y\right>$ and we claim that $M=W$.

Let $\varrho:W\rightarrow\sym Z$ denote the natural projection map in
Lemma~\ref{lem}.
If either $n$ or $m$ is odd then
$\left<(1,2),(1,2,\ldots,m)\right>\leq M\varrho$; otherwise
$\left<(1,2),(2,3,\ldots,m)\right>\leq M\varrho$. As
$$
\left<(1,2),(1,2,\ldots,m)\right>=\left<(1,2),(2,3,\ldots,m)\right>=\sym Z,
$$
we obtain that $M\varrho=\sym Z$. As $S=\ker\varrho$,
in order to prove that $M=W$, it suffices to show that $S\leq M$.

Next we claim that $A\leq M$. If $n=2$, then  $A=1$, and so we may assume that
$n\geq 3$. First we suppose
that $n=3$. In this case $A\cong(C_3)^m$ and so $A$ can be
viewed as a $W$-module over $\F_3$, and, in particular, it can be viewed as a $\sym
Z$-module over the same field. In fact, $A$ is the natural permutation module for $\sym Z$.
Since $M\varrho=\sym Z$, the intersection $A\cap M$ is an
$\sym Z$-submodule of $A$.
Now $y^2=((1,2,3),(1,2,3),\id,\ldots,\id)$ and so $y^2\in A\cap M$, but,
if $m\geq 3$, then
$y^2$ is not an element of either of the two proper submodules listed in Lemma~\ref{replem}.
Therefore $A\cap M=A$, and so $A\leq M$ follows when $n=3$ and $m\geq 3$. The
case $(n,m)=(3,2)$ can be checked using a computer algebra package such as
{\sf GAP}~\cite{gap} or {\sc Magma}~\cite{Magma}.

Let us assume that $n=4$, that is, $Y=\{1,\ldots,4\}$,
and, as above, we may also assume without loss of
generality that $m\geq 3$. Note
that $\alt Y$ admits the decomposition $\alt Y=U\rtimes V$ where
$U=\left<(1,2)(3,4),(1,3)(2,4)\right>$ and $V=\left<(1,2,3)\right>$.
Further,
$U$ can be considered as an irreducible $V$-module over $\F_2$.

Let  $m$ be odd and set $z_1=(x^my^2)^2$.
As
$$
x^my^2=((1,3,4),(1,3,4),(1,2),\ldots,(1,2)),
$$
we obtain that
$z_1=(x^my^2)^2=((1,4,3),(1,4,3),\id,\ldots,\id)$.
As $M\varrho=\sym Z$,  the subgroup $M$  has an element of the form
$\pi(2,3)$ where $\pi\in S$. Then set $z=z_1^{\pi(2,3)}$ and compute that
$z=z_1^{\pi(2,3)}=(\sigma_1,\id,\sigma_3,\id,\ldots,\id)$ where $\sigma_1$ and
$\sigma_3$ are three-cycles in $\sym Y$.
Now set
$$
w_1=(x^m)^{y^4}x^m=((1,2)(3,4),(1,2)(3,4),\id,\ldots,\id)
$$
and consider the element
$w=w_1^{z}w_1$. Now, as the first component of $z\in S$ is a non-trivial three
cycle, the element $w$ is of the form $w=(\sigma,\id,\ldots,\id)$ where
$\sigma\in\{(1,2)(3,4),(1,3)(4,2),(1,4)(2,3)\}$.

Assume now that $m$ is even and set $z=(x^{m-1}y^2)^4$.
Easy computation yields that
$$
z=(x^{m-1}y^2)^4=(\id,(1,3,4),\id,\ldots,\id)$$
and that
$$
w=(x^{m-1})^{y^4}x^{m-1}=
(\id,(1,2)(3,4),\id,\ldots,\id).
$$

As $V$ is irreducible on $U$,
the computation above shows that the  $\left<z\right>$-submodule
generated by $w$ coincides with $U\times
1\times\cdots\times 1$ of $m$ is odd, and $1\times U\times
1\times\cdots\times 1$ if $m$ is even. Thus $U\times
1\times\cdots\times 1\leq M$ in the former case, and $1\times U\times
1\times\cdots\times 1\leq M$ in the latter.
As $M\varrho=\sym Z$, we also obtain that
$U^m\leq M$.

Now the quotient $A/U^m\cong V^m\cong (C_3)^m$  can be considered as a permutation module for $\sym
Z$ over $\F_3$, and as $(M\cap A)/U^m\unlhd M/U^m$ and $M\varrho=\sym Z$, we
obtain that $(M\cap A)/U^m$ is a $\sym Z$-submodule. However, the image
of the element
$z$ above is not in either of the proper submodules listed in
Lemma~\ref{replem}, and so we obtain that $(M\cap A)/U^m=A/U^m$ which shows
that $A\leq M$.

Hence we have shown that $A\leq M$ if $n\leq 4$. Assume now that $n\geq 5$.
In this case $\alt Y$ is a non-abelian
finite simple group and so $A$ is a non-abelian
characteristically simple group.

Set $Q=M\cap S$. Clearly, $Q\leq S$ and
$Q\unlhd M$. Let us show that $Q$
is a subdirect subgroup of $S=\sym Y^m$.
If $n$ or $m$ is odd then set
$z=x^m=((1,2),(1,2),\ldots,(1,2))$; otherwise
set $z=x^{m-1}=(\id,(1,2),(1,2),\ldots,(1,2))$. Further,
$$
y^2=((1,2,\ldots,n),(1,2,\ldots,n),\id,\ldots,\id).
$$
 Therefore $z,\
y^2\in Q$. Let $\varrho_2$ denote the second coordinate
projection $\varrho_2:S\rightarrow\sym Y$. Since $z\varrho_2=(1,2)\in
Q\varrho_2$ and $y^2\varrho_2=(1,2,\ldots,n)\in
Q\varrho_2$, we obtain that $Q\varrho_2=\sym
Y$. Simple computation shows that if $\tau\in Q$ and
$\pi\sigma\in M$ then
$\tau^{\pi\sigma}\varrho_{2\sigma}=\tau^\pi\varrho_2$.
Therefore
$$
\sym Y=Q\varrho_2=Q^{\pi\sigma}\varrho_{2\sigma}=Q\varrho_{2\sigma}.
$$
As $Q\varrho_2=\sym Y$, it follows that $Q\varrho_i=\sym Y$ for all
$i$, and so $Q$ is a subdirect subgroup of $S=\sym Y^m$.
Consider the commutator subgroup $Q'$. Since $S'=A$, we obtain that $Q'\leq
A$. Further, as $Q\varrho_i=\sym Y$, we also obtain, for all $i$, that
$Q'\varrho_i=\alt Y$. Therefore $Q'$ is a subdirect subgroup of $A=\alt
Y^m$. Set $R=M\cap A$. Since $Q'\leq R$ and $Q'$ is a subdirect subgroup, we
obtain that so is $R$; that is, by Lemma~\ref{scott1},
$R$ is a direct product of disjoint strips.

Let $S$ be a strip in $R$ and let $\mathcal S\subseteq\{1,\ldots,m\}$ be the
support of $S$. We claim that $\mathcal S$ is a block for the action of $\sym
Z$. Indeed, if $\sigma\in\sym Z$ then there is some $\pi\in S$ such
that $\pi\sigma\in M$. Then $S^{\pi\sigma}$ is strip in $R$ and
so either $S=S^{\pi\sigma}$ or $S\cap S^{\pi\sigma}=1$. The support of
$S^{\pi\sigma}$ is $\mathcal S^{\sigma}$. Thus either $\mathcal
S^{\sigma}=\mathcal S$ or $\mathcal S^\sigma\cap
\mathcal S=\emptyset$, which shows that $\mathcal S$ is a block, as required. Since $\sym Z$ is
primitive on $\{1,\ldots,m\}$ either $|\mathcal S|=1$ or $\mathcal S=\{1,\ldots,
m\}$. If the
latter holds, then
$Q$ is a strip. This, however, is impossible. Indeed,
if $m=2$ and $n$ is odd, then set $z=xy=(\id,(1,2)(1,2,\ldots,n))$; if
$m=2$ and $n$ is even then set $z=xy^2=((1,2,\ldots,n),(1,2)(1,2,\ldots,n))$; if
 $m\geq 3$, then set
$z=y^2=((1,2,\ldots,n),(1,2,\ldots,n),\id,\ldots,\id)$.
Then in all cases $z^2\in R$,  but
$z^2$ is not in a full strip (if $m=2$ and $n$ is even then note that the first
component of $z^2$ is of order $n/2$ and the second is of order $n-1$). Thus
$|\mathcal S|=1$, and so $A\leq M$.

This completes the proof of the claim that that $A\leq M$ for all $n$ and $m$.

Note that $A$ is a normal subgroup of $W$ and let $x\mapsto \widehat x$
denote the natural homomorphism $W\rightarrow W/A$. If $H\leq W$, then
$\widehat H$ denotes the image $HA/A$ of $H$.
Then $\widehat W\cong C_2\wr\sym Z$ and $\widehat S\cong (C_2)^m$.
We claim that $\widehat S\leq \widehat M$.
Note that $\widehat M\cap\widehat S\unlhd \widehat M$.
If
$\pi\sigma\in\widehat M$ and $u\in \widehat M\cap\widehat S$ then $u^{\pi\sigma}
=u^\sigma$. As
$u^{\pi\sigma}\in \widehat M\cap\widehat S$ we obtain that $u^\sigma\in
\widehat M\cap\widehat S$ which shows that $\widehat M\cap\widehat S$ is a
$\sym Z$-submodule of $\widehat S\cong (C_2)^m$.
It is clear that $\widehat S$ is the natural permutation module for $\sym Z$
over $\F_2$. Lemma~\ref{replem} lists the non-trivial proper submodules $U_1$
and $U_2$ of
$\widehat S$.

If both $n$ and $m$ are even then $\widehat
x^{m-1}=(0,1,\ldots,1)$ and this
element is not in either $U_1$ or $U_2$. Hence $\widehat S\leq \widehat M$.
If this is not the case, then $\widehat x^m=(1,1,\ldots,1)$ which shows that
$U_1\leq \widehat M$.
If $n$ is even and $m$ is odd then $\widehat y^2=(1,1,0,\ldots,0)$ and so
$U_2\leq \widehat M$. Therefore in this case $U_1\oplus U_2\leq \widehat M$
follows.

Suppose that $n$ is odd.
If $m=2$ then
$xy=(\id,(1,2)(1,2,\ldots,n))$. Thus $\widehat{xy}$ is not in
either of the proper submodules listed in Lemma~\ref{replem}. Hence
$\widehat S\leq\widehat M$ follows in this case.
If $m\geq 3$, then
$$
xy=
(\id,(1,2),\id,\ldots,\id,(1,2,\ldots,n))(2,3\ldots,m).
$$
Let $\pi=(\id,(1,2),\id,\ldots,\id,(1,2,\ldots,n))$ and
$\sigma=(2,3,\ldots,m)$ so that $xy=\pi\sigma$.
As $\sigma^{m-1}=1$, it follows from Lemma~\ref{lem} that
$$
(xy)^{m-1}=\pi\pi^{\sigma^{-1}}\ldots\pi^{\sigma^{-m+2}}.
$$
Thus $(xy)^{m-1}$ is of the form $(\id,\pi_2,\ldots,\pi_m)$ where,
for $i=2,\ldots,m$, the permutation $\pi_i$ is either $(1,2)(1,2,\ldots,n)$ or
 $(1,2,\ldots,n)(1,2)$, that is, $\pi_i$ is a cycle with length $n-1$. As $n$
is odd, $\pi_i\not\in\alt Y$, and so
$(\widehat {xy})^{m-1}=(0,1,1,\ldots,1)$. Now if $m$ is even then
$(\widehat{xy})^{m-1}$ is not an element of $U_1$ or $U_2$, and so $\widehat
S\leq \widehat M$
follows also in this case. If $m$ is odd then this shows that $U_2\leq
\widehat M$, and
as we proved above that $U_1\leq \widehat M$, it follows that $\widehat S\leq \widehat
M$.

Thus we have shown that $S\leq M$ as required.
As explained above, $M=W$ now follows.
\end{proof}

The main result of the paper can now be proved.

\begin{proof}[The proof of Theorem~$\ref{main}$]
As usual, we assume without loss of generality that
 $X=Y\times Z$ where $Y=\{1,\ldots,n\}$, $Z=\{1,\ldots,m\}$
and $\P$ is the partition
$$
\{\{(1,1),\ldots,(n,1)\},\ldots,\{(1,m),\ldots,(n,m)\}\}.
$$
By Lemmas~\ref{wreq},~\ref{l1},~\ref{l2} and Theorem~\ref{wrth},
\begin{multline*}
\rank{\transs X\P}=\rank{(\transs X\P:\sym{X,\P})}+\rank{\sym{X,\P}}\\=
\rank{(\trans Y\wr\trans Z:\sym Y\wr\sym Z)}+\rank{\sym Y\wr\sym Z}=4.
\end{multline*}
The assertions concerning $\rank\Gamma(X,\P)$ and $\rank\Sigma(X,\P)$ can be proved very
similarly.
\end{proof}

\section*{Acknowledgment}
We acknowledge with gratitude some conversations with John D. Dixon, Peter M. Neumann and Joseph J. Rotman.

The first author was partially supported by FCT and FEDER, Project
POCTI-ISFL-1-143
of Centro de Algebra da Universidade de Lisboa, and by FCT and PIDDAC
through the project PTDC/MAT/69514/2006.

The second author is grateful to the Centro de \'Algebra da Universidade de
Lisboa for the invitation and the hospitality; he was also supported by the Hungarian Scientific
Research Fund (OTKA)
grant F049040.

\end{document}